\documentclass[10pt]{article}

\usepackage[T2A]{fontenc}
\usepackage[cp1251]{inputenc}
\usepackage[russian, english]{babel}
\usepackage{algpseudocode}

\newcommand{\sech}{{\rm sech}\,}
\newcommand{\csch}{{\rm csch}\, }
\newcommand{\eid}{{\bf EID}}
\newenvironment{trm}[1]{ \vskip 0.3cm \par {\bf Theorem #1} \it}{\vskip 0.3cm}
\newenvironment{proof}{\par $\bullet$ }{$\bullet$ \vskip 0.3cm}

\title{On Euler-Imshenetsky-Darboux transformation of second-order linear differential
  equations}
\author{L.M. Berkovich, S.A. Evlakhov, \\ Samara State University \\
        443011, Samara, Acad. Pavlov str.,\,1\\
        E-mail: berk@ssu.samara.ru}
\date{}

\begin{document}

\maketitle

\begin{abstract}
	It is shown, how to generate infinite sequences of differential equations of the second order 
  based on some standard equations, using Euler-Imshenetsky-Darboux (\eid) transformation.  
  For all this, fac\-to\-ri\-za\-tions of differential operators and operational identities are used.
  Some ge\-ne\-ra\-li\-za\-ti\-ons of integrable cases of Schr\"odinger's equations are finded.
  The example of integrable equation with liouvillian coefficients, that, apparently, can not be
  solved by Singer and Kovacic's algorithm (and its modifications), was built. The alogorithm
  and the program for solving conctructed classes of equations, was realized in REDUCE system.
  Corresponding procedure GENERATE is addendum to the ODESOLVE procedure of REDUCE system.
  The results of using GENERATE procedure (REDUCE 3.8) was compared with results  of DSolve
  procedure (Maple 10).
  Though, algorithm, based on \eid\ transformation is not an alternative to the unviersal
  ones, in the borders of its applicability, its very powerful.
\end{abstract}

\section{Introduction}

This papaer is dedicated to the classical Euler-Imshenetsky-Darboux transformation [1,2,3]
and its applications to the solving of differential equation.
First of all, lets make some general remarks

Let's consider incomplete equations
$$
	y'' + a_{0}y = 0, \quad a_{0}(x) \in C(I),\quad  I=(a,b). \eqno(1.1)
$$
Any complete equation
$$
	y'' +a_{1}y' + a_{0}y = 0, \quad a_{1} \in C^{1}(I),\quad a_{0}\in C(I) \eqno (1.1')
$$
can be reduced to (1.1) using the substitution
$$
	y=\exp(-\frac{1}{2}\int a_{1} dx)Y.
$$
Indeed, in this case we shall get incomplete equation 
$$
	Y'' + A_{0}(x)Y=0, \quad A_{0}=a_{0} -\frac{1}{4}a_{1}^{2} -\frac{1}{2}a_{1}',
$$
so, equations of the form (1.1) will be considered hereafter. Euler-Imshenetsky-Darboux
transformation (\eid) is  transformation of the form\footnote{Transformation (1.2) is well known as Darboux transformation and also B\"acklund
transformation. We propose to name it \eid, that is tribute to the historicas justice (see [4]).}
$$
  z= \beta(x)y' + \alpha(x)y, \quad \beta(x),\ \alpha(x) \in C^{2}(I). \eqno(1.2)
$$
\eid\ transformation is one of the two most important transformations of the linear equations.
Another important transformation is Kummer-Liouville's transformation, that has a form
$$
	y=v(x)z, \quad dt=u(x)dx,
$$
where $v(x)$ and $u(x)$ --- are sufficiently smooth functions 
(you can read about this transformation in detail in [4]).

{\bf Statement of the problem:} it requires to bring (1.1) to the form 
$$
	z'' + b_{0}(x)z=0, \quad b_{0}(x) \in C(I), \eqno(1.3)
$$ 	
using reversible transformation (1.2).
In other wods,
\begin{itemize}
	\item by the giving equations (1.1) and (1.3) find transformation (1.2). 
\end{itemize}
This problem allows another equivalent  formulatings:
\begin{itemize}
	\item by the given (1.1) and (1.2) find (1.3);
	\item by the given (1.2) and (1.3) find (1.1).
\end{itemize}

\subsection{Statement of the problem in the matrix form}
	Let's pass from the scalar equations (1.1)--(1.3) to the matrix ones. Introduce designations
$$
	Y=\left(
	\begin{array}{c}
		y \\
		y'\\
	\end{array} \right), \quad
	Z=\left(
	\begin{array}{c}
		z \\
		z' \\
	\end{array}\right), \quad
	A = \left(
	\begin{array}{cc}
		0 & 1 \\
		-a_{0} & 0 \\
	\end{array}\right), \quad
	B = \left(
	\begin{array}{cc}
		0 & 1 \\
	 -b_{0} & 0 \\
	\end{array}\right). \eqno(1.4)
$$

Instead of (1.1)--(1.3) we shall get:
$$
	Y'=AY, \eqno(1.5)
$$
$$
	Z=TY, \quad 
	T = \left(
	\begin{array}{cc}
			\alpha & \beta \\
			\alpha' - \beta a_{0} & \alpha + \beta' \\
	\end{array} \right), \quad 
	\det T \not = 0, \eqno(1.6)
$$
$$
	Z' =BZ. \eqno(1.7)
$$
It is required to find the transformation(1.6), that reducing (1.5) to (1.7).
\begin{trm}{1}{\rm [5, 6](see  also [4]).}
	Next propositions are equivalent: 
	\begin{description}
		\item[a)] Equation (1.1) can be reduced to (1.3) by transformation \eid\ (1.2).
		\item[b)] System (1.5) can be reduced to  (1.7) by transformation (1.6).
		\item[c)] Matrices $A,B$ and $T$ in formulas (1.4), (1.6) are such, that the next condition is held
			$$
				T' = BT - TA. \eqno(1.8)
			$$
		\item[d)] Functions $\beta(x)$ and $\alpha(x)$ satisfies the system of the equations
			$$
				\alpha'' + (b_{0} - a_{0})\alpha - 2a_{0}\beta' - \beta a_{0}' = 0, \eqno(1.9)
			$$
			$$
				\beta'' + (b_{0} - a_{0})\beta + 2\alpha' = 0. \eqno(1.10)
			$$
		\item[e)] The system of the equations (1.9), (1.10) lets first integral (FI)
			$$
				\alpha\beta' -\beta\alpha' + \alpha^{2}+a_{0}\beta^{2} = K,  \eqno(1.11)
			$$
			\begin{center}
				or (in equivalent form)
			\end{center}
			$$
				\alpha\beta' + \beta\alpha' + a^{2} + \beta\beta'' + b_{0}\beta^{2} = K, \eqno(1.12)
			$$
			where $K$ --- is arbitrary constant, $K \not = 0$: otherwise
			\eid\ transformation (see formula (1.6)) will be degenerated.
		\item[f)] Equation (1.3) by the transformation
			$$
				y=\frac{\alpha + \beta'}{K}z - \frac{\beta}{K}z' \eqno(1.13)
			$$
			can be reduced to (1.1).
		\item[g)] System (1.7) by the transformation
			$$
				Y=T^{-1}Z, \quad T^{-1} = \left(
				\begin{array}{cc}
					\frac{\alpha + \beta'}{K} & - \frac{\beta}{K} \\
              & \\
					\frac{a_{0}\beta - \alpha'}{K} & \frac{\alpha}{K} \\
				\end{array} \right) \eqno(1.14)
			$$
			can be reduced to (1.5).
		\item[h)] Between differential operators, that characterize direct and reverse \eid\ conversion 
		consequently
			$$
				P=\beta D + \alpha, \quad Q=\beta D - \alpha - \beta ', \eqno(1.15)
			$$
			there are next commutative  relations
			$$
				(QP+K)Q = Q(PQ+K), \eqno(1.16)
			$$
			$$
				P(QP+K) = (PQ+K)P, \eqno(1.17)
			$$
			where $K$  is the first integral of the form  (1.11), or (1.12).
	\end{description}
\end{trm}

Theorem 1 connects  corresponding results of the works [5,6]. Besides, in [6] 
irreversible transformation (1.2) was also considered.
The following proof is taken from [4].

\begin{proof}
	a) $\Rightarrow$ b). Immidiate corollary, because of the designations (1.4).

	b) $\Rightarrow$ c). Substituting (1.6) in (1.7), we shall get consecutively 
	$T'Y~+~TY'=BTY,~ T'Y + TAY =BTY$, from which come to (1.8).

	c) $\Rightarrow$ d). Let's write equation (1.8) in expanded form
	$$
		T' = \left(
		\begin{array}{cc}
			\alpha' & \beta' \\
			\alpha'' - \beta' a_{0} - \beta a_{0}' & \alpha' + \beta'' \\
		\end{array} \right) = 
	$$
	$$
		=\left(
		\begin{array}{cc}
			0 & 1 \\
			-a_{0} & 0 \\
		\end{array}
		\right)
		\left(
		\begin{array}{cc}
			\alpha & \beta \\
			\alpha' - \beta a_{0} & \alpha + \beta' \\
		\end{array}
		\right) -
		\left(
		\begin{array}{cc}
			\alpha & \beta \\
			\alpha' - \beta a_{0} & \alpha + \beta' \\
		\end{array}
		\right)
		\left(
		\begin{array}{cc}
			0 & 1 \\
			-b_{0} & 0 \\
		\end{array}
		\right),		
	$$
	from which follows the system of the equation (1.9), (1.10).

	d) $\Rightarrow$ e). Multiply (1.9) by $\beta$,  (1.10) by $\alpha$ and subtract first equation
  from second: $\alpha\beta'' - \beta\alpha'' + 2\alpha\alpha' + 2a_{0}\beta\beta' + a_{0}'\beta = 0$,
	from which follows FI of the form (1.11). Multiplying (1.10) by $\beta$ and adding with (1.11), 
  we shall get (1.12).

	e) $\Rightarrow$ f). From (1.2) by virtue of (1.1) the next equation follows
	$$
		z' = (\beta' + \alpha)y' + (\alpha' - \beta a_{0})y,
	$$
	(which coefficients corresponds to the lower row of the matrix $T$).
		
	Solving the system, consisting of (1.2) and  finded equation, relative to $y$, in 
  virtue of  (1.11) we shall come to (1.13).

	f) $\Rightarrow$ g). Reversing the transformation(1.6), we shall come to (1.14)
  in virtue of the relation
	$$
		\det T = K;
	$$	

	g) $\Rightarrow$ h). Indeed, 
	$\beta^{2}(D^{2} + a_{0}) = QP + K,\ \beta^{2}(D^{2} + b_{0})= PQ + K$, where $K=$FI. 
  But in this case	commutational correlations  (1.16), (1.17) are evident.

	h) $\Rightarrow$ a). Condition (1.15) include not only direct, but also reverse
   transformation \eid .
\end{proof}

\subsection{On solutions of linear equations  \\of second order}

Linear differential equation of the second order (1.1'),
that allows factorization of the form 
$(D - \alpha_{2})(D - \alpha_{1})y=0$, has general solution 
$$
  y(x)=e^{\int \alpha_{1}dx}(c_{1} + c_{2}\int e^{\int (\alpha_{2} - \alpha_{1})dx}dx).
$$
In particlular, incomplete equation  (1.1) allows the factorization \\
$(D + \alpha)(D -\alpha)y=0$ and has general solution
$$
   y(x) = e^{\int \alpha dx}(c_{1} + c_{2}\int e^{-2\int \alpha dx}dx)
$$
(more detail about this see in [4]).

\section{Procedure of <<generating>> \\ of the differential equations \\ using \eid\ transfomation}

Later on, instead of (1.2) we shall consider transformation
$$
	z = y' - \alpha y \quad (\beta=1, \ \alpha \rightarrow -\alpha). \eqno(2.1)
$$

From (1.9), (1.10) in virtue of (2.1) follows relations
$$
	\alpha'' + (b_{0} - a_{0})\alpha + a_{0}' = 0, \quad b_{0}=a_{0} + 2\alpha'. \eqno(2.2)
$$
Equation (1.11) will have a form
$$
	\alpha' + \alpha^{2}+ a_{0} =K \eqno(2.3).
$$
Though, (2.1) is the special case of the transformation (1.2), it is of interest as
 in the examination of the spectral problems both in generation of the
 remarkable sequence of the related equations, generated by the equations (1.1).

\begin{trm}{2}{\rm (Euler-Imshenetsky-Darboux)\ [4]}. 
  Equation (1.1), using transformation (2.1) generates the following sequence of equations
  $$
		y_{n}'' + a_{n}y_{n} = 0, \eqno(2.4)
	$$
	$$
		a_{n} = a_{0} + 2 \sum_{s=1}^{n} \alpha_{s-1}', \eqno(2.5)
	$$
	where  $\alpha_{s-1}$ satisfies the Riccati's equation
	$$
		\alpha_{s-1}' + \alpha_{s-1}^{2} + a_{s-1} = \lambda_{s-1}.	
	$$
	If denote
	$$
		\alpha_{s-1} = (\ln \tilde y_{s-1})' = \frac{\tilde y_{s-1}'}{\tilde y_{s-1}},
	$$
	where $\tilde y_{s-1}$ is eigenfunction of the equation
	$$
		y_{s-1}'' + (a_{s-1} - \lambda)y_{s-1} = 0,
	$$
	that corresponding  to the eigenvalue $\lambda=\lambda_{s-1}$, than $a_{n}$ will have a form 
	$$
		a_{n}=a_{0} + 2\sum_{s=1}^{n}\left( \frac{\tilde y_{s-1}'}{\tilde y_{s-1}}\right)',
	$$
	and solution  $y_{n}(x)$ of the equation (2.3) can be represented in one of the
  following three forms:
	$$
		y_{n} = \prod_{s=n}^{1}(D - \alpha_{s-1})y_{0} = 
    \prod_{s=n}^{1}(D - \frac{\tilde y_{s-1}'}{\tilde y_{s-1}})y_{0}=
		\tilde y_{n-1}\prod_{s=n}^{1}\left(D \frac{\tilde y_{s-2}}{\tilde y_{s-1}}\right) y_{0},
   \eqno(2.6)
$$
where $y=y_{0},\
 \alpha_{0} = \alpha,\ \tilde y_{-1} = 1$.
\end{trm}
We observe that partially contents of the theorem 3 can be found in  [1--3]. The
method of the proof is taken from [4].
\begin{proof}
	Starting equation  (1.1) by the transformation
	$$
		y_{1} = (D-\alpha_{0})y_{0} = \left(D - \frac{\tilde y_{0}'}{\tilde y_{0}}\right)y_{0},
	$$
	where $\alpha_{0}(x)$ satisfies the equation
	$$
		\alpha' + \alpha^{2} + a_{0} = \lambda
	$$
	when $\lambda = \lambda_{0}$ and $\tilde y_{0}$ is eigenfunction of the equation
	$$
		y'' + (a_{0} - \lambda)y=0, \eqno(2.7)
	$$
	that corresponding to the eigenvalue $\lambda_{0}$, is reducing to the form 
	$$
		y_{1}'' + \left(a_{0} + 2\left(\frac{\tilde y_{0}'}{\tilde y_{0}}\right)'\right)y_{1} = 0. \eqno(2.8)
	$$
	Next, by the transformation
	$$
 		y_{2} = \left(D - \frac{\tilde y_{1}'}{\tilde y_{1}}\right)y_{1},
	$$
	where $\tilde y_{1}$ is eigenfunction of the equation
	$$
		y_{1}'' + \left(a_{0} + 2\left( \frac{\tilde y_{0}'}{\tilde y_{0}} \right)' - \lambda \right) y_{1} = 0,
	$$
	that corresponding to the eigenvalue $\lambda_{1}$, (2.7) is reducing to the equation
	$$
		y_{2}'' + \left(a_{0} + 2\left(\frac{\tilde y_{0}'}{\tilde y_{0}}\right)' + 2
			\left(\frac{\tilde y_{1}'}{\tilde y_{1}} \right)' \right)y_{2} = 0.
	$$
	Proceeding this, we shall find, that the substitution
	$$
		y_{n} = \left(D - \frac{\tilde y_{n-1}'}{\tilde y_{n-1}}\right)y_{n-1}, \eqno(2.9)
	$$
	where $\tilde y_{n-1}$ is eigenfunction of the equation
	$$
		y_{n-1}'' + \left(a_{0} + 2\sum_{k=0}^{n-2}\left(\frac{\tilde y_{k}'}{\tilde y_{k}}\right)' 
			- \lambda \right)y_{n-1} = 0,
	$$
	corresponding to the eigenvalues $\lambda_{n-1}$, reduces to the equation
	$$
		y_{n}'' + \left(a_{0} + 2\sum_{k=0}^{n-1}\left(\frac{\tilde y_{k}'}{\tilde y_{k}}\right)'
			\right)y_{n} = 0,
	$$
	that corresponding to (2.3). Substitution (2.9) reduces to the first of the equality (2.5).
  The next representation of the factorization is used in proof of the last equality 
  in formula (2.6) (see [4, $\S$ 5.5]):
	$$
		\prod_{s=n}^{1}(D - \alpha_{s-1}) = \exp(\int \alpha_{n-1}dx)\prod_{s=n}^{1}[D\exp (-\int(\alpha_{s-1}
		-\alpha_{s-2})dx)], \quad \alpha_{-1}=0
	$$
\end{proof}

Theorem 2, in fact, gives us an algorithm for generating the sequences of the multiplying
equations, that can be realized in the interactive mode. The algorithm's feature is that it
simultaneously construct  both equation and form of its solution.

\subsection{Example}
We consider the equation
$$
	y_{0}'' + \lambda^{2}y_{0}=0. \eqno(2.10)
$$
When $\lambda_{0}=0$ let $\tilde y_{0}=x$. Then $\tilde y_{0}'=1,\ \tilde y_{0}'/\tilde y_{0}=1/x,$
$(\tilde y_{0}'\tilde y_{0})'=-1/x^2$. We shall get an equation 
$$
	y_{1}'' + (\lambda^{2} - \frac{2}{x^2})y_{1} =0, \eqno(2.11)
$$
 it's general solution has a form
$$
  y_{1}(x) = \left(D - \frac{1}{x}\right)(c_{1}\cos \lambda x + c_{2}\sin \lambda x). \eqno(2.12)
$$
Accept (2.11) (received from  (2.12) when $c_{1}=0, c_{2}=1, \lambda=1$) function
$$
	\tilde y_{1}(x) = (D - \frac{1}{x})\sin x = \cos x - \frac{1}{x}\sin x,  
$$
as the partial solution.
Then
$$
	\frac{\tilde y_{1}'}{\tilde y_{1}} = -\frac{1}{x} - \frac{x\sin x}{x\cos x  - \sin x},\quad
	\left(\frac{\tilde y_{1}'}{\tilde y_{1}}\right)' = \frac{1}{x^2} - \frac{x^{2} - \sin^{2} x}{(x\cos x  - \sin x)^{2}}.
$$
So, in the result, we shall get an equation
$$
	y_{2}'' + (\lambda^{2} - \frac{2 (x^{2} - \sin^{2} x)}{(x\cos x - \sin x)^{2}})y_{2}=0. \eqno(2.13)
$$
General solution of the (2.11) can be represented in one of the next forms:
$$
  y_{2}(x) = (D + \frac{1}{x} + \frac{x\sin x}{x\cos x - \sin x})(D - \frac{1}{x})(c_{1}\cos \lambda x + c_{2}\sin \lambda x),
$$
or
$$
	y_{2}(x) = (D^{2} + \frac{x \sin x}{x\cos x - \sin x}D - \frac{\sin x}{x\cos x - \sin x})(c_{1}\cos \lambda x + c_{2}\sin \lambda x).
$$
Applying REDUCE\footnote{We have used REDUCE 3.8.}, we shall get final expression 
for $y_{2}(x)$.
$$
	y_{2}(x)=c_{1}\frac{-\lambda^{2} x\cos \lambda x \cos x + (\lambda^{2} -1)\cos \lambda x\sin x - \lambda x \sin \lambda x \sin x}{x\cos x - \sin x} +
$$
$$
 +c_{2}\frac{- \lambda^2 x \cos x \sin \lambda x   +
  (\lambda^2 - 1)\sin \lambda x \sin x + \lambda x \cos \lambda x \sin x}{x\cos x - \sin x}
$$
 Authors can't solve this equation, using Singer and Kovacic's algorithm [8, 9].
 DSolve procedure (Maple 10) also can't solve this equation.

Using theorem  2, one can enlarge the number of examples, but we shall be interested in
those of the sequences, whose common terms can be counctructed, using the identities, that 
are given below.

\section{Special class \\of the related equations}

\begin{trm}{3} If common term of the equation sequence, generated by the EID transformation
has a form
$$
	y_{n}'' + \left[a_{0} + n(n+1)\left(\frac{\tilde y_{0}'}{\tilde y_{0}}\right)'\right]y_{n} = 0,
		\quad \frac{\tilde y_{0}'}{\tilde y_{0}} = \alpha_{0}, \eqno(3.1)
$$
then its solution  (general or partial) can be represented in the form
$$
	y_{n}=\prod_{k=n}^{1}(D-k\frac{\tilde y_{0}'}{\tilde y_{0}})y_{0}=\tilde y_{0}^{n}
	(D\frac{1}{\tilde y_{0}})^{n}y_{0}, \eqno(3.2)
$$
or
$$
	y_{n} = \prod_{k=n}^{0}\left(D - k \frac{\tilde y_{0}'}{\tilde y_{0}} \right)D^{-1}y_{0} =
		\tilde y_{0}^{n+1}(\tilde y_{0}^{-1}D)^{n+1}D^{-1}y_{0}, \eqno(3.3)
$$
where $y_{0}$ is common or partial solution of the equation (1.1), and $\tilde y_{0}$  
is partial solution of the equation (2.6), corresponding to the eigenvalue $\lambda_{0}$.
\end{trm}
\begin{proof}
Let's accept, that $\tilde y_{k} = \tilde y_{0}^{k+1}.$ Then
$$
	a_{n} = a_{0} + 2\left[(\ln \tilde y_{0})''+ (\ln \tilde y_{1})'' + \ldots + (\ln \tilde y_{n-1})''
		 \right] = a_{0} + n(n+1)\left(\frac{\tilde y_{0}'}{\tilde y_{0}}\right)'.
$$
And
$$
	y_{n} = \left(D - n \frac{\tilde y_{0}'}{\tilde y_{0}} \right) \ldots 
	\left(D - 2\frac{\tilde y_{0}'}{\tilde y_{0}}\right)\left(D-\frac{\tilde y_{0}'}{\tilde y_{0}}\right)y_{0}.
$$
This transformation can be rewrited in the another form
$$
	y_{n} = \prod_{k=n}^{0}\left(D - k \frac{\tilde y_{0}'}{\tilde y_{0}} \right)D^{-1}y_{0} =
		\tilde y_{0}^{n+1}(\tilde y_{0}^{-1}D)^{n+1}D^{-1}y_{0}.
$$
Indeed,  $y_{n}$ in factorized form can be transformated:
$$
	\prod_{k=n}^{1}(D-k\frac{\tilde y_{0}'}{\tilde y_{0}}) =\tilde y_{0}^{n}
	(D\frac{1}{\tilde y_{0}})^{n}.
$$
Let's multiply this identity by  $D$ in right
$$
	\prod_{k=n}^{0}(D - k \frac{\tilde y_{0}'}{\tilde y_{0}}) = \tilde y_{0}^{n} (D \frac{1}{\tilde y_{0}})^{n}D =
	\tilde y_{0}^{n+1}(\frac{1}{\tilde y_{0}} D)(\frac{1}{\tilde y_{0}} D)^{n}=
	\tilde y_{0}^{n+1}(\frac{1}{\tilde y_{0}} D)^{n+1}
$$
from where we shall have (3.2), (3.3).
\end{proof}
Let's find members of the class (3.1),
  where $\tilde y_{0}$ satisfies the equation (2.7) where $\lambda = \lambda_{0}$. 
  Because factorization
	$D^{2} + a_{0} - \lambda_{0} = (D + \alpha_{0})(D - \alpha_{0})$ takes plase, we shall have
	$a_{0} = \lambda_{0} - \alpha_{0}^{2} - \alpha_{0}'$. Correspondingly 
	$a_{k} = \lambda_{k} - \alpha_{k}^{2} - \alpha_{k}'$. And also
$$
	a_{k} = a_{0} + 2\sum_{s=1}^{k}\alpha_{s-1}' = \lambda_{0} - \alpha_{0}^{2} - \alpha_{0}'
		+ 2\sum_{s=1}^{k}s\alpha_{0}'.
$$
We also consider, that $\alpha_{s-1} = s\alpha_{0}$. After simple calculations
$$
	a_{k} = \lambda_{k} - (k+1)^{2}\alpha_{0}^{2} - (k+1)\alpha_{0}'.
$$
As a result we shall have an equation
$k(k+2)(\alpha_{0}' + \alpha_{0}^{2}) = \lambda_{k} - \lambda_{0}$.
When $\lambda_{k} = \lambda_{0}$ we can find $\alpha_{0}' + \alpha_{0}^{2} = 0$.
 This equation has a solution $\alpha_{0} = 1/(x+c)$, $c$ --- is an arbitrary constant. Than
$a_{0}(x) = 0, \lambda_{k} = \lambda_{0} = 0$, i.e. (2.7) will have a form $y'' = 0$. Accepting
$\lambda_{k} =(k+1)^{2}$, we shall find solutions of the equation $\alpha_{0}'+\alpha_{0}^{2} = 1$,
 that  has a form
$\alpha_{0} = \tanh(x+c),\ \alpha_{0}=\coth(x+c)$. Accepting $\lambda_{k} = - (k+1)^{2}$, we shall
find for equation  $\alpha_{0}' + \alpha_{0}^{2} = -1$ solutions of the form $\alpha_{0} = - \tan(x+c), 
\alpha_{0} = \cot(x+c)$. In accordance with theorem 3, we shall construct some important examples.

\section{Operational identities and some types \\ of the Schr\"odinger's equations}

Let's consider equations, generated by equations

$$
	y_{0}'' - \lambda y_{0} = 0. \eqno(4.1)
$$

$1^{0}$. $\lambda_{0} = 0,\ \tilde y_{0} = x,\ \lambda_{k}=0,\ \tilde y_{k} =x^{k+1},\ 
	a_{n} = -\lambda - \frac{n(n+1)}{x^{2}}$, i.e. we shall have an equation
$$
  y_{n}'' - \left(\lambda + \frac{n(n+1)}{x^{2}}\right)y_{n} = 0. \eqno(4.2)
$$
By virtue of operational identity
$$
  \prod_{k=n}^{1}\left(D - \frac{k}{x}\right)=x^{n}\left(D\frac{1}{x}\right)^{n},
$$
that connects factorization of the differential operators with iteration of the
same differential operator of the first order, the next identity is also correctly
$$
	\prod_{k=n}^{0}(D - \frac{k}{x})=x^{n+1}(\frac{1}{x}D)^{n+1}.
$$
One can get it, multiplying previous identity by $D$ in right. Indeed, one can go to
the required identity, transposing factors in the right part of the identity and using
the associative rule. So, general (partial) solution of the equation (4.2) can be express
in the form:
 
$$
	y_{n} = x^{n+1}\left(\frac{1}{x}D\right)^{n+1}y_{0}, \eqno(4.3)
$$
where $y_{0}$  is general (partial) solution of the equation (4.1).

General solution of the equation (4.1) and also its partial solutions 
 (with concrete values of $c_{1}$, $c_{2}$) 
can be express with exponents, hyperbolic and trigonometric functions:
$$
	y_{0}=c_{1}\exp(\sqrt{\lambda}x) + c_{2}\exp(-\sqrt{\lambda}x), \eqno(4.4)
$$
$$
	y_{0}=c_{1}\cosh(\sqrt{\lambda}x) + c_{2}\sinh(\sqrt{\lambda}x),\quad \lambda>0, \eqno(4.5)
$$
$$
	y_{0}=c_{1}\cos(\sqrt{-\lambda}x) + c_{2}\sin(\sqrt{-\lambda}x),\quad \lambda<0, \eqno(4.6)
$$
where $c_{1}$ and $c_{2}$  are arbitrary constants. 

When $\lambda=0$ differential operator  (DO), corresponding to (4.2), allows
factorizations
$$
	D^{2} - \frac{n(n+1)}{x^{2}} = \left(D + \frac{n+1}{x}\right)\left(D - \frac{n+1}{x}\right)=
	\left(D - \frac{n}{x}\right)\left(D + \frac{n}{x}\right).
$$
And corresponding equation
$$
  y'' -\frac{n(n+1)}{x^2}y=0 \eqno(4.7)
$$ has general solution
$$
  y=c_{1}x^{n+1} + c_{2}x^{-n} \eqno(4.8)
$$ (see $\S$ 1.2). However,  (4.7) is Euler's equation, and its solution  (4.8) 
can be finded, using known Euler's substitution $x=e^{t}$.

$2^{0}$. $\lambda_{0}=1,\ \tilde y_{0}=\cosh x,\ \lambda_{k} =(k+1)^{2},\ \tilde y_{k} =\cosh^{k+1}x,\
a_{n}=-\lambda - \frac{n(n+1)}{\cosh^{2}x}$, i.e. we shall get an equation
$$
	y_{n}'' - \left(\lambda - \frac{n(n+1)}{\cosh^{2}x}\right)y_{n} = 0. \eqno(4.9)
$$
Let's multiply operational identity
$$
 \prod_{k=n}^{1}\left(D - k \tanh x\right) = \cosh^{n}x\left(D \frac{1}{\cosh x} \right)^{n}
$$
by $D$ in right. Then transpose factors in the right part of the identity, using associative
rule. We shall get an identity
$$
 \prod_{k=n}^{0}\left(D - k \tanh x\right) = \cosh^{n+1}x\left(\frac{1}{\cosh x} D\right)^{n+1}.
$$
(4.9) has a solution
$$
	y_{n} = \cosh^{n+1} x \left(\frac{1}{\cosh x}D\right)^{n+1}y_{0},\eqno(4.10)
$$
where $y_{0}$ can be taken in accordance with formulas (4.4)--(4.6).
When $\lambda=(n+1)^{2}$ DO, corresponding to (4.9), allows a factorization
$$
	D^{2} - (n+1)^{2} + \frac{n(n+1)}{\cosh^{2}x} = (D + (n+1)\tanh x)(D -(n+1)\tanh x),
$$
and differential equation
$$
	y'' -(n+1)^{2}y + \frac{n(n+1)}{\cosh^{2}x}y =0 \eqno(4.11)
$$
has general solution
$$
	y = \cosh^{n+1} x\left(c_{1} + c_{2}\int \frac{dx}{\cosh^{2(n+1)}x}\right) \eqno(4.12)
$$
(see p. 1.2), where integral in formula (4.12)
can  be expressed in elementary functions (see  [7], $\S$  1.2):
$$
\int \frac{dx}{\cosh^{2(n+1)}x}=
$$
$$=\frac{\sinh x}{2n+1}\left\{
\sech^{2n+1}x +\sum_{k=1}^n\frac{2^kn(n-1)\ldots(n-k+1)}{
(2n-1)(2n-3)\ldots(2n-2k+1)}\sech^{2n-2k+1}x\right\}.
$$

$3^{0}$. $\lambda_{0} = 1,\ \tilde y_{0} = \sinh x,\ \lambda_{k}=(k+1)^{2},\ \tilde y_{k} = \sinh^{k+1}x,
a_{n}=-\lambda - \frac{n(n+1)}{\sinh^{2}x}$, i.e. we shall have an equation
$$
	y_{n}'' - \left(\lambda + \frac{n(n+1)}{\sinh^{2}x} \right)y_{n} = 0. \eqno(4.13)
$$
Let's multiply operational identity 
$$
 \prod_{k=n}^{1}\left(D - k \coth x\right) = \sinh^{n}x\left(D \frac{1}{\sinh x} \right)^{n}
$$
by $D$ in right and transpose factors in the right part of the identity, using the associative rule. We shall have an identity
$$
 \prod_{k=n}^{0}\left(D - k \coth x\right) = \sinh^{n+1}x\left(\frac{1}{\sinh x} D\right)^{n+1}.
$$
(4.13) has a solution
$$
 y_{n} = \sinh^{n+1} x \left(\frac{1}{\sinh x}D\right)^{n+1}y_{0}, \eqno(4.14)
$$
where $y_{0}$ is one of the values (4.4)--(4.6).
When $\lambda=(n+1)^{2}$ DO, corresponding to (4.13), allows factorization
$$
	D^{2} - (n+1)^{2} - \frac{n(n+1)}{\sinh^{2}x} = (D + (n+1)\coth x)(D - (n+1)\coth x),
$$
and corresponding differential equation 
$$
	y'' -(n+1)^{2}y - \frac{n(n+1)}{\sinh^{2}x}y =0 \eqno(4.15)
$$
has general solution
$$
	y = \sinh^{n+1} x\left(c_{1} + c_{2}\int \frac{dx}{\sinh^{2(n+1)}x}\right) \eqno(4.16)
$$
(see $\S$ 1.2), where in accordance with [7, p. 110]
$$
\int\frac{dx}{\sinh^{2(n+1)}x}=
$$
$$
=\frac{\cosh x}{2n+1}\left\{-\csch^{2n+1}x +\sum_{k=1}^n(-1)^{k-1}\frac{2^kn(n-1)\ldots(n-k+1)}{(2n-1)(2n-3)\ldots(2n-2k+1)}\csch^{2n-2k+1}
x\right\}.
$$

$4^{0}$. $\lambda_{0} = -1,\ \tilde y_{0} = \cos x,\ \lambda_{k} = - (k+1)^{2},\ 
\tilde y_{k}=\cos^{k+1}x,\ a_{n} = - \lambda - \frac{n(n+1)}{\cos^{2}x}$, i.e. 
we shall have an equation
$$
	y_{n}'' - \left(\lambda + \frac{n(n+1)}{\cos^{2} x} \right)y_{n} = 0. \eqno(4.17)
$$
Let's multiply operational identity 
$$
 \prod_{k=n}^{1}\left(D - k \tan x\right) = \cos^{n}x\left(D \frac{1}{\cos x} \right)^{n}
$$
by $D$ in right and transpose factors in the right part of identity, using the associative rule.
We shall have an identity:
 $$
 \prod_{k=n}^{0}\left(D - k \tan x\right) = \cos^{n+1}x\left(\frac{1}{\cos x} D\right)^{n+1}.
$$
(4.17) has a solution
$$
 y_{n} = \cos^{n+1} x \left(\frac{1}{\cos x}D\right)^{n+1}y_{0},\eqno(4.18)
$$
where $y_{0}$ is one of the values (4.4)--(4.6)
When $\lambda=-(n+1)^{2}$ DO, corresponding to (4.17), allows factorization
$$
	D^{2} + (n+1)^{2} - \frac{n(n+1)}{\cos^{2}x} = (D - (n+1)\tan x)(D + (n+1)\tan x),
$$
and corresponding differential equation 
$$
	y'' + (n+1)^{2}y - \frac{n(n+1)}{\cos^{2}x}y =0 \eqno(4.19)
$$
has general solution
$$
	y = \cos ^{n+1} x\left(c_{1} + c_{2}\int \frac{dx}{\cos^{2(n+1)}x}\right). \eqno(4.20)
$$
(see $\S$ 1.2), where in accordance with  [7, p. 149]
$$
\int \frac{dx}{\cos^{2(n+1)}x}=
$$
$$=\frac{\sin x}{2n+1}\left\{\sec^{2n+1}x
 +\sum_{k=1}^n\frac{2^kn(n-1)\ldots (n-k+1)}{(2n-1)(2n-3)\ldots(2n-2k+1)}\sec^{2n-2k+1}x\right\}.
$$

$5^{0}$. $\lambda_{0} = -1,\ \tilde y_{0} = \sin x,\ \lambda_{k} = - (k+1)^{2},\ 
\tilde y_{k} = \sin^{k+1} x,\ a_{n} = -\lambda - \frac{n(n+1)}{\sin^{2}x}$, 
i.e. we shall have an equation:
$$
	y_{n}'' - \left(\lambda + \frac{n(n+1)}{\sin^{2}x}\right)y_{n} = 0. \eqno(4.21)
$$
Let's multiply operational identity
$$
 \prod_{k=n}^{1}\left(D - k \cot x\right) = \sin^{n}x\left(D \frac{1}{\sin x} \right)^{n}
$$
by $D$ in right and transpose factors in the right part of the identity,
using associative rule. We shall have an identity
$$
 \prod_{k=n}^{0}\left(D - k \cot x\right) = \sin^{n+1}x\left(\frac{1}{\sin x} D\right)^{n+1}.
$$
Equation (4.21) has general solution
$$
 y_{n} = \sin^{n+1} x \left(\frac{1}{\sin x}D\right)^{n+1}y_{0}, \eqno(4.22)
$$
where $y_{0}$ is one of the values (4.4)--(4.6).
When $\lambda = - (n+1)^{2}$ DO, corresponding to the (4.21), allows factorization
$$
	D^{2} + (n+1)^{2} - \frac{n(n+1)}{\sin^{2} x} = (D + (n+1)\cot x)(D - (n+1)\cot x),
$$
and corresponding differential equation
$$
	y'' + (n+1)^{2}y - \frac{n(n+1)}{\sin^{2}x}y =0 \eqno(4.23)
$$
has general solution
$$
	y = \sin ^{n+1} x\left(c_{1} + c_{2}\int \frac{dx}{\sin^{2(n+1)}x}\right) \eqno(4.24)
$$
(see $\S$ 1.2), where in accordance with [7, p. 148]
$$
\int\frac{dx}{\sin^{2(n+1)}x}=
$$
$$=-\frac{\cos x}{2n+1}\left\{\csc^{2n+1}x
 +\sum_{k=1}^n\frac{2^kn(n-1)\ldots(n-k+1)}{(2n-1)(2n-3)\ldots(2n-2k+1)}\csc^{2n-2k+1}x\right\}.
$$

{\bf Remark.} Equation (2.7) is known as Sturm-Liouville's  equation, and also
linear Schr\"odinger's equation.  Equations  (4.2), (4.9), (4.13), (4.17) and (4.21) 
are known as Schr\"odinger's equations with soluble (integrable) potentials, 
but we have constructed them by the another way. When $n=1$ this equations are
partial cases of the Lame's equation,
$$
	y'' - (2\wp(x) + \nu)y=D,
$$
where  $\wp(x)$  is the elliptical function of Weirstrass . This partial cases 
corresponds to the elementary values of Weirstrass  function. At that parameter $\lambda$ 
for (4.9), (4.13), (4.17) and (4.21) when $n=1$ differs from $\nu$  in Lame's equation.
In [4, p.67--70] this was explored, using another methods.

\section{Generalized operational identities\\ and corresponding Schr\"odinger equations \\}

5.1 Equation
$$
	y_{n}'' - \left(\lambda + \frac{n(n+1)a^{2}}{(ax+b)^{2}}\right)y_{n}=0 \eqno(5.1)
$$
by virtue of identity
$$
	\prod_{k=n}^{0}\left(D - k\frac{a}{ax+b}\right)=(ax+b)^{n+1}\left(\frac{1}{ax+b}D\right)^{n+1} \eqno(5.2)
$$
has general solution
$$
	y_{n}=(ax+b)^{n+1}\left(\frac{1}{ax+b}D\right)^{n+1}y_{0}, \eqno(5.3)
$$
where $y_{0}$ is one of the values  (4.4)--(4.6).
When $\lambda=0$ we shall have an equation
$$
	y_{n}'' - \frac{n(n+1)a^{2}}{(ax+b)^{2}}y_{n} = 0, \eqno(5.4)
$$
that allows factorizations
$$
	D^{2} - \frac{n(n+1)a^{2}}{(ax+b)^{2}} = \left(D - \frac{na}{ax+b}\right)
	\left(D + \frac{na}{ax+b}\right) =
$$
$$=\left[D + \frac{(n+1)a}{ax+b}\right]\left[D - \frac{(n+1)a}{ax+b}\right]
$$
and has general solution
$$
	y=c_{1}(ax+b)^{n+1} + c_{2}(ax+b)^{-n}. \eqno(5.5)
$$

5.2. Let
$$
	\lambda_{0}=m^{2},\quad \tilde y_{0} = ae^{mx} + be^{-mx}, \quad \lambda_{k} = m^{2}(k+1)^{2},
	\quad \tilde y_{k} = (ae^{mx}+be^{-mx})^{k+1}.
$$
Then
$$
	\frac{\tilde y_{0}'}{\tilde y_{0}} =
\frac{m(ae^{mx} - be^{-mx})}{ae^{mx}+be^{-mx}}, \quad
\left(\frac{\tilde y_0'}{\tilde y_0}\right)'=
	 \frac{4abm^{2}}{(ae^{mx} + be^{-mx})^{2}}.
$$
Equation
$$
	y_{n}'' - \left(\lambda - \frac{4abm^{2}(n+1)n}{(ae^{mx} + be^{-mx})^{2}}\right)y_{n} =0, \quad \lambda \not = 0 \eqno(5.6)
$$
by virtue of identity
$$
	\prod_{k=n}^{0}\left(D - km\frac{ae^{mx} - be^{-mx}}{ae^{mx} + be^{-mx}}\right) =
	 \left(ae^{mx} + be^{-mx}\right)^{n+1}\left(\frac{1}{ae^{mx} + be^{-mx}}D \right)^{n+1} \eqno(5.7)
$$
has general solution
$$
	y_{n}=\left(ae^{mx} + be^{-mx}\right)^{n+1}\left(\frac{1}{ae^{mx}+be^{-mx}}D\right)^{n+1} y_{0}, \eqno(5.8)
$$
where $y_{0}$ is one of the values (4.4)--(4.6).

Let $\lambda =m^{2}(n+1)^{2}$. By virtue of
$$
	D^{2} - m^{2}(n+1)^{2} + \frac{4abm^{2}n(n+1)}{(ae^{mx} + be^{-mx})^{2}}=
$$
$$
	=\left[D + m(n+1)\frac{ae^{mx} - be^{-mx}}{ae^{mx} + be^{-mx}}\right]
	\left[D - m(n+1)\frac{ae^{mx} - be^{-mx}}{ae^{mx} + be^{-mx}}\right]
$$
equation
$$
	y'' - m^{2}(n+1)^{2}y + \frac{4abm^{2}n(n+1)}{(ae^{mx} + be^{-mx})^{2}}y=0 \eqno(5.9)
$$
has general solution
$$
	y = \left(ae^{mx} + be^{-mx}\right)^{n+1}\left[c_{1} + c_{2}\int\left(ae^{mx} + be^{-mx}\right)^{-2(n+1)} dx\right]. \eqno(5.10)
$$

Let
$$
	\lambda_{0} = m^{2}; \quad \tilde y_{0} = a\cosh mx + b \sinh mx, \quad \lambda_{k} = m^{2}(k+1)^{2},
$$
$$
	\tilde y_{k} = \left(a\cosh mx + b\sinh mx \right)^{k+1}.
$$
Then
$$
	\frac{\tilde y_{0}'}{\tilde y_{0}} = m \frac{a \sinh mx + b \cosh mx}{a\cosh mx + b\sinh mx},\quad
	\left(\frac{\tilde y_{0}'}{\tilde y_{0}}\right)' = m^{2} 
	\frac{a^{2} - b^{2}}{(a \cosh mx + b\sinh mx)^{2}}.
$$

5.3. Equation
$$
  y_{n}'' - \left(\ell - \frac{n(n+1)m^{2}\left(a^{2} - b^{2}\right)}{\left(a\cosh mx + b\sinh mx\right)^{2}}\right)y_{n}=0 \eqno(5.11)
$$
by the virtue of the operational identity
$$
	\prod_{k=n}^{0}\left(D - km\frac{a\sinh mx + b\cosh mx}{a \cosh mx + b\sinh mx}\right) =
$$
$$
	= \left(a\cosh mx + b\sinh mx\right)^{n+1}\left(\frac{1}{a\cosh mx + b\sinh mx}D\right)^{n+1} \eqno(5.12)
$$
has general solution
$$
	y_{n} = \left(a\cosh mx + b\sinh mx\right)^{n+1}\left(\frac{1}{a\cosh mx + b\sinh mx}D\right)^{n+1}
	y_{0}, \eqno(5.13)
$$
where $y_{0}$ is one of the values (4.4)--(4.6).

Let $	\lambda = m^{2}(n+1)^{2}$. By the virtue of the operational identity
$$
	D^2 - m^{2}(n+1)^{2} + \frac{n(n+1)m^{2}(a^{2} - b^{2})}{(a\cosh mx + b\sinh mx)^{2}} =
$$
$$
	=\left[D + m(n+1)\frac{a\sinh mx + b\cosh mx}{a\cosh mx + b \sinh mx}\right]
	\left[D - m(n+1)\frac{a\sinh mx + b\cosh mx}{a\cosh mx + b \sinh mx}\right]
$$
Equation
$$
	y'' - m^{2}(n+1)^{2}y + \frac{n(n+1)m^{2}(a^{2} - b^{2})}{(a\cosh mx + b\sinh mx)^{2}}y =0 \eqno(5.14)
$$
has general solution
$$
	y=(a\cosh mx + b \sinh mx)^{n+1}\left[c_{1} + c_{2}\int\left(a\cosh mx + b\sinh mx\right)^{-2(n+1)} dx \right]. \eqno(5.15)
$$

Let
$$	
  \lambda_{0}=-m^{2}
$$
$$
	\tilde y_{0} = a \cos mx + b\sin mx, \quad \lambda_{k} = - m^{2}(k+1)^{2}, \quad
		\tilde y_{k} = (a \cos mx + b \sin mx)^{k+1}.
$$
Then
$$
	\frac{\tilde y_{0}'}{\tilde y_{0}}=m\frac{a\sin mx + b \cos mx}{a \cos mx + b \sin mx}, \quad
	\left( \frac{\tilde y_{0}'}{\tilde y_{0}}\right)' = - m^{2} \frac{a^{2} + b^{2}}{(a\cos mx + b\sin mx)^{2}}.
$$

5.4. Equation
$$
	y_{n}'' - \left(\lambda + \frac{m^{2}(a^{2} + b^{2})n(n+1)}{(a\cos mx + b\sin mx)^{2}}\right)y_{n}=0 \eqno(5.16)
$$
by the virtue of the operational identity
$$
	\prod_{k=n}^{0}\left[D - km\frac{-a\sin mx + b\cos mx}{a \cos mx + b\sin mx}\right] =
$$
$$
	=(a \cos mx + b \sin mx)^{n+1}\left(\frac{1}{a\cos mx + b \sin mx}D\right)^{n+1} \eqno(5.17)
$$
has general soltuion
$$
	y_{n} = (a\cos mx + b \sin mx)^{n+1}\left(\frac{1}{a\cos mx + b \sin mx}D\right)^{n+1}
	y_{0}, \eqno(5.18)
$$
where $y_{0}$ is one of the values (4.4)---(4.6).

Let $	\lambda=-m^{2}(n+1)^{2}.$ By the virtue of the factorization
$$
	D^2 + m^{2}(n+1)^{2} - \frac{n(n+1)m^{2}(a^{2} + b^{2})}{(a\cos mx + b\sin mx)^{2}} =
$$
$$
	=\left[D + m(n+1)\frac{-a\sin mx + b\cos mx}{a\cos mx + b \sin mx}\right]
	\left[D - m(n+1)\frac{-a\sin mx + b\cos mx}{a\cos mx + b \sin mx}\right]
$$
Equation
$$
	y'' + m^{2}(n+1)^{2}y - \frac{n(n+1)m^{2}(a^{2} + b^{2})}{(a\cos mx + b\sin mx)^{2}}y =0 \eqno(5.19)
$$
has general solution
$$
	y=(a\cos mx + b \sin mx)^{n+1}\left[c_{1} + c_{2}\int\left(a\cos mx + b\sin mx\right)^{-2(n+1)} dx \right]. \eqno(5.20)
$$

\section{Algorithm for solving \\ of the constructed equations} 

We consider following equations:

\noindent equation with rational coefficients of the form:
$$
	y_{n}'' - \left(\ell + \frac{n(n+1)a^{2}}{(ax+b)^{2}}\right)y_{n}=0;
$$
\noindent equation with exponential coefficients of the form:
$$
	y_{n}'' - \left(\ell - \frac{4abm^{2}(n+1)n}{(ae^{mx} + be^{-mx})^{2}}\right)y_{n} =0, \quad \ell \not = 0;
$$
\noindent equation with hyperbolic coefficients of the form:
$$
  y_{n}'' - \left(\ell - \frac{n(n+1)m^{2}\left(a^{2} - b^{2}\right)}{\left(a\cosh mx + b\sinh mx\right)^{2}}\right)y_{n}=0,\quad \ell \not = 0;
$$
\noindent equation with trigonometric coefficients of the form:
$$
	y_{n}'' - \left(\ell + \frac{m^{2}(a^{2} + b^{2})n(n+1)}{(a\cos mx + b\sin mx)^{2}}\right)y_{n}=0, \quad\ell \not = 0.
$$

Generating equation has a form:
$$
	y_{0}'' - \ell y_{0} = 0.
$$
Exponential type of the generating equations:
$$
  y_{0}=c_{1}\exp(\sqrt{l}x) + c_{2}\exp(-\sqrt{l}x).
$$
Hyperbolic type of the generating equation:
$$
  y_{0}=c_{1}\cosh(\sqrt{l}x) + c_{2}\sinh(\sqrt{l}x), \quad l>0.
$$
Trigonometric type of the generating  equation:
$$
  y_{0}=c_{1}\cos(\sqrt{-l}x) + c_{2}\sin(\sqrt{-l}x), \quad l<0.
$$

We use the type of the equation (lin, expon, hyp, trig 
for equations with rational, exponential, hyperbolic ant trigonometrical
coefficients correspndingly,), number n  of the equations in the sequence, parameters $a, b, m, l$,
the form of the generating equation's solution (expon, hyp and trig
for exponential, hyperbolic and trigonometric types of solution correspondingly)
and arbitrary constants $c_{1}$ и $c_{2}$.

Input: equation's form and it's solution.

\begin{center}
	\begin{tabular}{|c|c|}
		\multicolumn{2}{c}{{\bf Variables}} \\ \hline
		 $tp$ & The type of the equation\\ \hline
     $tpy_0$ & Form of the generating equaion's solution \\hline
     $n$ & The number of the equations in the sequence \\ \hline
     $a, b, m, l$ & Parameters (look above) \\ \hline
     $c_{1}, c_{2}$ & Arbitrary constants \\ \hline 
	\end{tabular}
\end{center}
\begin{algorithmic}[1]
	\State $y_{p}: =0$
	\If {$tpy_0 = {\rm expon}$} 
	   \State $y_n:=c_{1}\exp({\sqrt{l}x}) + c_{2}\exp({-\sqrt{l}x})$
  \ElsIf {$tpy_0 = {\rm trig}$} 
	  \State $y_n:=c_{1}\cos({\sqrt{-l}x}) + c_{2}\sin({\sqrt{-l}x})$
  \ElsIf {$tpy_0 = {\rm hyp}$} 
	  \State $y_n:=c_{1}\cosh({\sqrt{l}x}) + c_{2}\sinh({\sqrt{l}x})$
	\EndIf
	\If {общий случай} 
		\If {$tp={\rm lin} \quad {\rm and} \quad l \not = 0$} 
		  \State $y_n:=(ax + b)^{n+1}\left(\displaystyle{\frac{1}{(ax+b)}}D\right)^{n+1}y_n$
		\ElsIf {$tp={\rm lin} \quad {\rm and} \quad l = 0 $}
		  \State $y_n:=(ax+b)^{n+1}c_{1} + c_{2}(ax+b)^{-n}$
		\ElsIf {$tp={\rm expon} \quad {\rm and} \quad l \not = m^2(n+1)^2 $}
			\State $y_n:=(a\exp(mx) - b\exp(-mx))^{n+1}\left(\displaystyle{\frac{1}{a\exp(mx) + b\exp(-mx)}}D\right)^{n+1}y_n$
		\ElsIf {$tp={\rm expon} \quad {\rm and} \quad l = m^2(n+1)^2$}
			\State $y_n:=(a\exp(mx) + b\exp(mx))^{n+1}(c_{1} + c_{2}\int(a\exp(mx) + b\exp(-mx))^{-2(n+1)}dx)$
		\ElsIf {$tp={\rm trig} \quad {\rm and} \quad l \not = -m^2(n+1)^2$}
			\State $y_n:=(a\cos(mx) + b\sin(mx))^{n+1}\left(\displaystyle{\frac{1}{a\cos(mx) + b\sin(mx)}}D\right)^{n+1}y_n$
		\ElsIf {$tp={\rm trig} \quad {\rm and} \quad l = -m^2(n+1)^2$}
			\State $y_n=(a \cos (mx) + b\sin (mx))^{n+1}(c_{1} + c_{2}\int(a \cos (mx) + b \sin (mx) )^{-2(n+1)}dx)$
		\ElsIf {$tp={\rm hyp} \quad {\rm and} \quad l \not = m^2(n+1)^2$}
			\State $y_n=(a\cosh(mx) + b\sinh(mx))^{n+1}\left(\displaystyle{\frac{1}{a\cosh(mx) + b\sinh(mx)}}D\right)^{n+1}y_n$
		\ElsIf {$tp={\rm hyp} \quad {\rm and} \quad l = m^2(n+1)^2$}
 			\State $y_n=(a\cosh(mx) + b\sinh(mx))^{n+1}(c_{1} + c_{2}\int (a\cosh(mx) + b\sinh(mx))^{-2(n+1)}dx)$
		\EndIf
	\Else
	  \If {$tp = {\rm lin}$}
		  \State $ty_{n}:=(ax+b)^{n+1}$
		  \State $a_{n}:=l + n(n+1)a^{2}/(ax+b)^2$
		  \State $\alpha:=a/(ax+b)$
	  \EndIf
	  \If {$tp = {\rm expon} $} 
		  \State $ty_{n}:=(a \exp{(mx)} + b \exp{(-mx)})^{n+1}$
		  \State $a_{n}:=l - 4abn(n+1)m^{2}/(a\exp{(mx)} + b\exp{(-mx)})^{2}$
		  \State $\alpha:=  4abm^{2}/(a\exp{(mx)} + b\exp{(-mx)})^{2}$
		  \If {$l=m^{2}(n+1)^2$}
		     \State $y_{p}:=(a\exp{(mx)} + b\exp{(-mx})^{n+1}(c_{1}+ c_{2}\int(a\exp{(mx)} + b\exp{(-mx)})^{-2(n+1)}dx)$
   		\EndIf
	  \EndIf
	  \If {$tp = {\rm trig} $}
		  \State $ty_{n}:=(a\cos(mx) + b\sin(mx))^{n+1}$
		  \State $a_{n}:=l + n(n+1)m^{2}(a^{2} + b^{2})/(a\cos(mx) + b\sin(mx))^{2}$
		  \State $\alpha:=m(-a\sin(mx) + b\cos(mx))/(a\cos(mx) + b\sin(mx))$
		   \If {$l=m^{2}(n+1)^{2}$}
		    \State  $y_{p}:=(a\cos(mx) + b\sin(mx))^{n+1}(c_{1} + c_{2}\int(a\cos(mx) + b\sin(mx))^{-2(n+1)}dx)$
		   \EndIf
	  \EndIf
  	\If {$tp = {\rm hyp} $}
	  	\State $ty_{n}:=(a\cosh(mx) + b\sinh(mx))^{n+1}$
		  \State $a_{n}:=l - n(n+1)m^{2}(a^{2} - b^{2})/(a\cosh(mx) + b\sinh(mx))^{2}$
		  \State $\alpha:=m(a\sinh(mx) + b\cosh(mx))/(a\cosh(mx) + b\sinh(mx))$
	  	 \If {$l=m^{2}(n+1)^{2}$}
	  	    \State  $y_{p}:=(a\cosh(mx) + b\sinh(mx))^{n+1}(c_{1} + c_{2}\int(a\cosh(mx) + b\sinh(-mx))^{-2(n+1)}dx)$
	  	 \EndIf
  	\EndIf
	  \For{$i:=1$ to $n$}
	  	\State $tmp_{1}:=i\alpha y_{k}$
	  	\State $tmp_{2}:= d y_{n}/dx$
	  	\State $y_{n}:=tmp_{2} - tmp_{1}$
  	\EndFor
  \EndIf
\end{algorithmic}

\section{Examples}

Now you can see examples of the differential equations and their solutions, received
with GENERATE procedure (REDUCE 3.8, [10,11]) and DSolve procedure (Maple 10, [12]).

{\it Example 1.} Let's consider the following equation:
$$
y'' \left( x \right) - \left( l-2\, \left( 
\cosh \left( x \right)  \right) ^{-2} \right) y \left( x \right) = 0.
$$

Solution, received with DSolve procedure:
$$
   y(x)=\frac{1}{\sqrt{\sinh{2x}} (\cosh{2x} + 1)^{1/4}}\times
$$
$$
	\times c_{1}(2\cosh{2x} - 2)^{3/4} \, \mbox{hyperheom}
  \left( \left[-\frac{1}{2}\sqrt{l}, \frac{1}{2}\sqrt{l}\right], 
  \left[{-\frac{1}{2}}\right], \frac{1}{2}\cosh{2x} + \frac{1}{2}\right) +
$$
$$
  + \frac{1}{\sqrt{\sinh{2x}}}\Bigg( c_{2}(2\cosh{2x} - 2)^{3/4}(2\cosh{2x} + 2)^{1/4}\times
$$
$$
  \times \,\mbox{hyperheom}\left(\left[\frac{3}{2} + \frac{1}{2}\sqrt{l}, \frac{3}{2}-
  \frac{1}{2}\sqrt{l}\right], \left[\frac{5}{2}\right], \frac{1}{2}\cosh{2x} + \frac{1}{2}\right)\times
$$
$$
 \times (\cosh{2x} + 1) \Bigg).
$$
 
As we can see, solution of the comparatively simple equation is expressed in the hyperheometric
terms.

And this is solution, received with GENERATE procedure, expressed in the explicit form
through the elementary functions:
$$
  y_{1}=\frac{e^{2\sqrt{l}x} c_{1}\left(\sqrt{l}\cosh(x) - \sinh(x)\right) -
  c_{2}\left(\sqrt{l}\cosh(x) + \sinh(x)\right)}{e^{\sqrt{l}x} \cosh(x)}
$$

{\it Example 2:}
$$ 
  y'' \left( x \right) - \left( l+2\,{\frac {{a}^
 {2}+{b}^{2}}{ \left( a\cos \left( x \right) +b\sin \left( x \right) 
  \right) ^{2}}} \right) y \left( x \right) =0.
$$
DSolve can't finds solution:
$$
  y(x) = \mbox{DESol}\left(\left\{ \left(\frac{d^{2}}{dx^{2}} Y(x)\right) + 
  \left(-l - \frac{2(a^{2} + b^{2})}{(a\cos{x} + b\sin{x})^2}\right)Y(x)\right\}, \{ Y(x)\}\right).
$$ 

GENERATE performs the solution in the elementary functions:
$$
  y_{1}= \frac{e^{2\sqrt{l}x} c_{1}\left(a\sqrt{l}\cos(x)-b\cos(x) +
   b\sqrt{l}\sin(x) + a\sin(x)\right)}{e^{\sqrt{l}x}\left(a\cos(x) + b\sin(x)\right)} +
$$
$$ +\frac{
   c_{2}\left(-a\sqrt{l}\cos(x) - b\cos(x) - b\sqrt{l}\sin(x) + a\sin(x)\right)}
  {e^{\sqrt{l}x}\left(a\cos(x) + b\sin(x)\right)}
$$

{\it Example 3:}
$$
  y'' -  \left(l + \frac{n(n+1)a^{2}}{(ax+b)^2}\right)y=0.
$$
Solution, received with DSolve:
$$
y \left( x \right) ={\it \_C1}\, \left( ax+b \right) ^{-n}{e^{-{\frac 
{\sqrt {l} \left( ax+b \right) }{a}}}}{\rm hypergeom} \left( [-n],[-2
\,n],2\,{\frac {\sqrt {l} \left( ax+b \right) }{a}} \right)+
$$
$$ +{\it \_C2
}\, \left( ax+b \right) ^{n+1}{e^{-{\frac {\sqrt {l} \left( ax+b
 \right) }{a}}}}{\rm hypergeom} \left( [n+1],[2+2\,n],2\,{\frac {
\sqrt {l} \left( ax+b \right) }{a}} \right). 
$$
Solution, received with GENERATE:
$$
	y_{n}=(ax+b)^{n}\left(\frac{1}{(ax+b)}D\right)^{n+1} \left(\frac{c_{1}e^{2\sqrt{l}x} + c_{2}}{e^{\sqrt{l}x}}\right).
$$

{\it Example 4}:
$$
	y'' - \left(l + \frac{6}{\sin^{2}x} \right)y=0.
$$
Solution, received with DSolve:
{\footnotesize
$$
y \left( x \right) ={\it \_C1}\left( \frac {  \left( -4\,\cos
 \left( 2\,x \right) -8+2\,l\cos \left( 2\,x \right) -2\,l \right) 
\sinh \left( \sqrt {l}\arcsin \left( 1/2\,\sqrt {2\,\cos \left( 2\,x
 \right) +2} \right)  \right)}{-1+\cos \left( 2\,x \right) }\right. +
$$

$$ 
+ \left. \frac{-3\,\sqrt {2\,\cos \left( 2\,x \right) +
2}\cosh \left( \sqrt {l}\arcsin \left( 1/2\,\sqrt {2\,\cos \left( 2\,x
 \right) +2} \right)  \right) \sqrt {-2\,\cos \left( 2\,x \right) +2}
\sqrt {l}  }{-1+\cos \left( 2\,x \right) }\right)+ 
$$

$$
+{\it \_C2}\left(\frac {
  3\, \left(  \left( i\sqrt {l}+1/3\,l-2/3 \right) \cos
 \left( 2\,x \right) +i\sqrt {l}-1/3\,l-4/3 \right)}{-1+\cos \left( 2\,x \right) }\right. \times
$$
$$ \times
 \frac{ \sqrt {-2\,\cos
 \left( 2\,x \right) +2}\cos \left(  \left( i\sqrt {l}-1 \right) 
\arcsin \left( 1/2\,\sqrt {2\,\cos \left( 2\,x \right) +2} \right) 
 \right)}{-1+\cos \left( 2\,x \right) }+ 
$$

$$
+ \frac{ -3\,\sqrt {2\,\cos \left( 2\,x \right) +2}\sin \left( 
 \left( i\sqrt {l}-1 \right) \arcsin \left( 1/2\,\sqrt {2\,\cos
 \left( 2\,x \right) +2} \right)  \right)}{-1+\cos \left( 2\,x \right)} \times
$$
$$
\times \left.\frac{\left(  \left( i\sqrt {l}+1
/3\,l-2/3 \right) \cos \left( 2\,x \right) -4/3-i\sqrt {l}-1/3\,l
 \right) }{-1+\cos \left( 2\,x \right)}\right).
$$}
Solution, received with GENERATE:
$$
 y_{2}=\frac{e^{2\sqrt{l}x}c_{1}\left(3\cos^2(x) - 
  3\sqrt{l}\cos(x)\sin(x) + l\sin^2(x) + \sin^2(x)\right)}
 {e^{\sqrt{l}x}\sin^2(x)} +
$$
$$
  + \frac{c_{2}\left(3\cos^{2}(x) + 3\sqrt{l}\cos(x)\sin(x)  + l\sin^2(x)
  + \sin^{2}(x)\right)}
 {e^{\sqrt{l}x}\sin^2(x)}
$$

\section{Conclusion}
Differential equations, that has theoretical or practical significance, can be 
divided into equivalence classes by the types of the used transformations. All equations,
that belongs to the same class are related and stay in relationship with standard equation,
that has gave them birth. Solutions of all equations, belongs to the same class, can
be expressed through the solutions of the sample equation. In the terms of
diffrential algebra (see [13]) it means, that if constructed equation has liouvillian
coefficietns, then it's Picar-Vessiot extension belongs to the differential field,
generated by the corresponding liouvillian coefficients. Sequential applying of the \eid\
transformation allows to find such integrable equations, that can't be solve
even in so powerful system as Maple, though it uses Singer and Kovacic's algorithm.
In our work Schr\"odinger's equations with soluble potentials, that generalaize well-known
Schr\"odinger's equations, depending of some parametrs, were constructed, using factorizations
of differential operators. Solutions (eigenfunctions) of this Schr\"odinger's equations are
expressed in terms of equations with constant coefficients. The constructed algorithm, that
finds solutions of the generated equations, also describes a process of their origin.
Though, it realized as GENERATE procedure in REDUCE 3.8, it can be realized in other
systems of the computer algebra. In the boundaries of its applicability, this procedure
is not less, than universal procedures.

\end{document}